\documentclass[12pt]{amsproc}
\usepackage{amsmath,amssymb,amscd}

\begin{document}

\title[Almost Normal Operators mod Hilbert--Schmidt]{Almost Normal Operators mod Hilbert--Schmidt \\ and the $K$-theory of the Algebras $E\Lambda(\Omega)$}
\author{Dan-Virgil Voiculescu}
\address{D.V. Voiculescu \\ Department of Mathematics \\ University of California at Berkeley \\ Berkeley, CA\ \ 94720-3840}
\thanks{Research supported in part by NSF Grant DMS 1001881.}
\date{}

\begin{abstract}
Is there a mod Hilbert-Schmidt analogue of the
BDF-theorem, with the Pincus $g$-function playing
the role of the index? We show that part of the
question is about the $K$-theory of certain Banach
algebras. These Banach algebras, related to Lipschitz
functions and Dirichlet algebras have nice Banach-space
duality properties. Moreover their corona algebras
are $C^*$-algebras.
\end{abstract}

\maketitle

\section{Introduction}
\label{sec1}

The BDF-theorem \cite{6} classifies, up to unitary equivalence, the normal elements of the Calkin algebra, by the spectrum and the index of the resolvent. If the ideal of compact operators is replaced by the trace-class, for operators with trace-class self-commutator, the Pincus $g$-function (\cite{7}, \cite{8}) is an $L^1$-function on ${\mathbb C}$ which extends the index of the essential resolvent.  The $g$-function has been related to algebraic $K$-theory by L.~G.~Brown (\cite{4}, \cite{5}) and in another direction, after work of J.~W.~Helton and R.~Howe (\cite{17}), the distribution to which the $g$-function gives rise, has been interpreted in terms of cyclic cohomology by A.~Connes (\cite{12}).

These developments around the $g$-function, were however not accompanied by a corresponding BDF-type result. In (\cite{27}, \cite{25}, \cite{26}) we formulated conjectures about operators with trace-class self-commutator, an affirmative answer to which would fill this gap. Besides the initial evidence in favor of these conjectures, there was no further progress. The situation is roughly that the $g$-function viewed in the cyclic cohomology framework covers the index part and our work on Hilbert--Schmidt perturbations of normal operators (\cite{24}) covers the part about trivial extensions, while the rest is wide open. The absence on the technical side of a normal dilation result which would correspond to the existence of inverses in Ext and which in the BDF context can be derived from the Choi--Effros completely positive lifting theorem, is a noted difficulty.

Our aim here is to decouple the normal dilation from the rest by introducing the algebras $E\Lambda(\Omega)$. In this way we are also able to bring $K$-theory to the study of this problem since we are led to the $K_0$-group of such an algebra.

The Banach $*$-algebras $E\Lambda(\Omega)$ are the natural framework to study operators with trace-class self-commutator which are obtained from compressions of normal operators to mod Hilbert--Schmidt reducing projections. Roughly $E\Lambda(\Omega)$, where $\Omega$ is a Borel subset of ${\mathbb C}$ is an algebra of operators in $L^2(\Omega,\lambda)$ with Hilbert--Schmidt commutators with the multiplication operators by Lipschitz functions, a construction reminiscent of Paschke-duality (\cite{21}).

The algebras $E\Lambda(\Omega)$ have nice properties as Banach algebras. They resemble the Lipschitz algebras of \cite{29}, up to the use of a Hilbert--Schmidt norm instead of a uniform norm, which is a feature of the Dirichlet algebras of non-commutative potential theory (\cite{1}, \cite{9}, \cite{10}). Actually the ideal ${\mathcal K}\Lambda(\Omega)$ of compact operators in $E\Lambda(\Omega)$ is a Dirichlet algebra and we show that $E\Lambda(\Omega)$ can be viewed both as the algebra of multipliers or as the bidual of ${\mathcal K}\Lambda(\Omega)$, when $\Omega$ is bounded. Since all this has the flavor of Banach algebra analogues of basic $C^*$-algebras, it is perhaps unexpected that the corona $E\Lambda(\Omega)/{\mathcal K}\Lambda(\Omega)$ which is the analogue of the Calkin algebra is really a $C^*$-algebra. Note, however, that while the Dirichlet algebra ${\mathcal K}\Lambda(\Omega)$ has the same simple $K$-theory as the algebra ${\mathcal K}({\mathcal H})$ of compact operators, the $K$-theory of $E\Lambda(\Omega)$ and hence of $E\Lambda(\Omega)/{\mathcal K}\Lambda(\Omega)$, which interests us in connection with operators with trace-class self-commutator, is certainly richer.

On the technical side an essential ingredient is the existence of a bounded approximate unit consisting of projections for ${\mathcal K}\Lambda(\Omega)$, which is a consequence of our work on norm-ideal perturbations of Hilbert-space operators (\cite{24}, \cite{28}).

Concerning the relation of the operator theory problems to the $K$-theory of the algebras $E\Lambda(\Omega)$, we should point out that while the $K$-theory problem is so to speak the operator theory problem minus the dilation problem, actually certain outcomes of the $K$-theory problem could provide a negative answer to the dilation problem. If the $K$-theory of $E\Lambda(\Omega)$ exhibits some integrality property making $K_0$ less rich this would answer in the negative the dilation problem.

In addition to the first section, which is the introduction, the paper has five more sections.

Section~2 contains background material about the conjectures about almost normal operators modulo Hilbert--Schmidt. Details of certain connections between these problems, left out previously, are included for the reader's convenience.

Section~3 introduces the algebras $E\Lambda(\Omega)$ and some of their basic properties. We also consider the ideal of compact operators ${\mathcal K}\Lambda(\Omega)$ of $E\Lambda(\Omega)$ and the Banach algebra $E\Lambda({\mathbb C})_0$ which is the inductive limit of the $E\Lambda(\Omega)$ for bounded sets $\Omega$.

In section~4 we look at the $K$-theory of the Banach algebras considered. We show that the problem about a mod Hilbert--Schmidt BDF-type theorem for almost normal operators is equivalent to the normal dilation problem plus the problem whether the $K_0$-group of $E\Lambda({\mathbb C})_0$ is isomorphic via the Pincus $g$-function to the group $L_{\mathop{\mathrm{re}}}^1({\mathbb C},\lambda)$ of real-valued $L^1$-functions with bounded support.

Section~5 returns to the algebras ${\mathcal K}\Lambda(\Omega)$, $E\Lambda(\Omega)$ and $(E/{\mathcal K})\Lambda(\Omega)$ and gives results about duality, multipliers and the relation to $C^*$-algebras.

Section~6 contains concluding remarks in several directions: the action of bi-Lipschitz homeomorphisms on the algebras, the problem about the center of $(E/{\mathcal K})\Lambda(\Omega)$, the relation to Dirichlet algebras and non-commutative potential theory, the possibility of similar constructions with other Schatten--von~Neumann classes ${\mathcal C}_p$ replacing the Hilbert--Schmidt class.

I would like to thank Jesse Peterson for useful comments about a possible connection to Dirichlet algebras.

\section{Background}
\label{sec2}

\noindent
{\bf 2.1.} If ${\mathcal H}$ is a separable infinite-dimensional Hilbert space over ${\mathbb C}$, then ${\mathcal B}({\mathcal H})$ will denote the bounded operators on ${\mathcal H}$ and ${\mathcal C}_p({\mathcal H})$ the {\em Schatten--von~Neumann $p$-class}. The $p$-norm $|\cdot|_p$ is $|T|_p = \mathop{\mathrm{Tr}}(T^*T)^{p/2}$. In particular, ${\mathcal C}_1({\mathcal H})$ is the {\em trace-class} and ${\mathcal C}_2({\mathcal H})$ is the {\em Hilbert--Schmidt class}.

\bigskip
\noindent
{\bf 2.2.} An operator $T \in B({\mathcal H})$ is {\em almost normal} if its {\em self-commutator} $[T^*,T]$ is in ${\mathcal C}_1({\mathcal H})$. Equivalently, if $T = A + iB$ with $A = A^*$, $B = B^*$ then $[A,B] \in {\mathcal C}_1$ since $2i[A,B] = [T^*,T]$. We shall denote by ${\mathcal A}{\mathcal N}({\mathcal H})$ the set of almost operators. Background material and references to the literature for many facts about operators with trace-class self-commutator can be found in the books \cite{11}, \cite{20}.

\bigskip
\noindent
{\bf 2.3.} If $T = A + iB \in {\mathcal A}{\mathcal N}({\mathcal H})$ and if $Q,R \in {\mathbb C}[X,Y]$ are polynomials in two commuting indeterminates, then since $\tilde{A},\tilde{B}$ the class of $A,B$ in ${\mathcal B}({\mathcal H})/{\mathcal C}_1({\mathcal H})$ commute, we shall also write $Q(A,B)$, $R(A,B)$ for elements in ${\mathcal B}({\mathcal H})$ so that $\widetilde{Q(A,B)} = Q(\tilde{A},\tilde{B})$, $\widetilde{R(A,B)} = R(\tilde{A},\tilde{B})$. Clearly these are only defined up to a ${\mathcal C}_1$ perturbation. The {\em Helton--Howe measure} $P_T$ of $T = A + iB \in {\mathcal A}{\mathcal N}({\mathcal H})$ (\cite{17}) is a compactly supported measure on ${\mathbb R}^2$ so that
\[
\mathop{\mathrm{Tr}}[Q(A,B),R(A,B)] = (2\pi i)^{-1} \int {\mathcal J}(Q,R)dP_T
\]
where
\[
{\mathcal J}(Q,R) = \frac {\partial Q}{\partial X} \frac {\partial R}{\partial Y} - \frac {\partial Q}{\partial Y} \frac {\partial R}{\partial X}\,.
\]
Then $\mathop{\mathrm{supp}} P_T \subset \sigma(T)$ and $P_T$ is absolutely continuous w.r.t.\ Lebesgue measure $\lambda$ and the Radon--Nikodym derivative $\frac {dP_T}{d\lambda} = g_T \in L^1({\mathbb R}^2)$ is the {\em Pincus principal function of $T$} (also called {\em Pincus $g$-function}).

\bigskip
\noindent
{\bf 2.4.} Let $R_1^+({\mathcal H}) = \{X \in B({\mathcal H}): X$ finite rank, $0 \le X \le 1\}$, which is a directed ordered set. Then the obstruction to the existence of quasicentral approximate units relative to the Hilbert--Schmidt class (\cite{24}) is
\[
k_2(T_1,\dots,T_n) = \liminf_{X \in R_1^+({\mathcal H})} \max_{1 \le j \le n} |[T_j,X]|_2.
\]
In \cite{27} we showed that: {\em if $T_1,T_2 \in {\mathcal A}{\mathcal N}({\mathcal H})$, $k_2(T_1) = 0$ and $T_1-T_2 \in {\mathcal C}_2$, then $P_{T_1} = P_{T_2}$ (or equivalently $g_{T_1} = g_{T_2}$ a.e.).}

\bigskip
\noindent
{\bf 2.5.} We recall two of the conjectures about almost normal operators (\cite{27} conjectures~3 and 4). Note that the second of these is a consequence of the first.

\bigskip
\noindent
{\bf Conjecture~3 in \cite{27}}. {\em If $T_1,T_2 \in {\mathcal A}{\mathcal N}({\mathcal H})$ are so that $P_{T_1} = P_{T_2}$ then there is a normal operator $N \in B({\mathcal H})$ and a unitary operator $U \in B({\mathcal H} \oplus {\mathcal H})$ so that $T_1 \oplus N - U(T_2 \oplus N)U^* \in {\mathcal C}_2$.}

\bigskip
If true, this statement would represent a kind of BDF-theorem with ${\mathcal A}{\mathcal N}({\mathcal H})$ and the Helton--Howe measure replacing the operators with compact self-commutator and respectively the index-data. Note also that the unitary equivalence is $\mathop{\mathrm{mod}} {\mathcal C}_2$ (not ${\mathcal C}_1$).

\bigskip
\noindent
{\bf Conjecture 4 in \cite{27}}. {\em If $T \in {\mathcal A}{\mathcal N}({\mathcal H})$ then there is $S \in {\mathcal A}{\mathcal N}({\mathcal H})$ and a normal operator $M \in B({\mathcal H} \oplus {\mathcal H})$ so that $T \oplus S - M \in {\mathcal C}_2$.}

\bigskip
This conjecture is an analogue of the existence of inverses in Ext in the analogue of the ``Ext is a group'' part of the BDF theorem. Note that the analogue of the results for trivial extensions (i.e., Weyl--von~Neumann theorem part) is covered by our results in \cite{24}. For the derivation of Conjecture~4 from Conjecture~3 one also uses the result of R.~V.~Carey and J.~D.~Pincus that every $L^1$-function is the $g$-function of some $T \in {\mathcal A}{\mathcal N}({\mathcal H})$.

\bigskip
\noindent
{\bf 2.6.} We would like to remark that Conjectures~3 and 4 in \cite{27} don't bring the essential spectrum of the almost normal operators into the discussion. With consideration of the essential spectrum $\sigma_e(T)$, one might ask if $P_{T_1} = P_{T_2}$ and $\sigma_e(T_1) = \sigma_e(T_2)$ would imply $T_1 - UT_2U^* \in {\mathcal C}_2$, form some unitary $U$, when $T_j \in {\mathcal A}{\mathcal N}({\mathcal H})$, $j = 1,2$.

We didn't discuss the possibility of such a strengthening, because it seems to have to do also with phenomena of another kind involving perturbations of isolated points in $\sigma(T)\backslash\sigma_e(T)$.

\bigskip
\noindent
{\bf 2.7.} A consequence of Conjecture~4 and hence also of Conjecture~3 is

\bigskip
\noindent
{\bf Conjecture 1 in \cite{27}}. {\em If $T \in {\mathcal A}{\mathcal N}({\mathcal H})$ then $k_2(T) = 0$.}

\bigskip
The proof which was omitted in \cite{27}, involves using {\em a result of} \cite{24}, {\em that $k_2(N) = 0$ for every normal operator $N$}. Indeed, if Conjecture~4 holds for $T$, then $T \in {\mathcal A}{\mathcal N}({\mathcal H})$ is unitarily equivalent $\mathop{\mathrm{mod}} {\mathcal C}_2$ to a compression $PN \mid P{\mathcal H}$ where $P = P^* = P^2$ is a projection, $N$ is normal and $[P,N] \in {\mathcal C}_2$. We infer that $k_2(T) = k_2(PN \mid P{\mathcal H})$. On the other hand $k_2(N) = 0$ implies there are $X_n \in R_1^+({\mathcal H})$, $X_n \uparrow I$ as $n \to \infty$, so that $\lim_{n \to \infty} |[X_n,N]|_2 = 0$. If $Y_n = PX_nP$ then $Y_n \in R_1^+$ and we have $Y_n \uparrow P$ as $n \to \infty$. We have
\[
|[Y_n,PNP]|_2 = |P[X_n,PNP]P|_2 \le |P[X_n,NP]P|_2 + |[I-X_n,[P,N]P]|_2.
\]
Since $[P,N]P \in {\mathcal C}_2$ and $I-X_n \downarrow 0$ we have $|[I-X_n,[P,N]P]|_2 \to 0$ as $n \to \infty$. On the other hand
\[
|P[X_n,NP]P|_2 \le |P[I-X_n,N]P|_2 + |P[N,P](I-X_n)P|_2
\]
which converges to $0$ as $n \to \infty$. Thus, Conjecture~1 holds for $T$, i.e., $k_2(T) = 0$.

\bigskip
\noindent
{\bf 2.8.} We will also need to recall some of the results for normal operators which follow from \cite{24}. Since $k_2(N) = 0$ for every normal operator $N$, we can use the kind of non-commutative Weyl--von~Neumann results in \cite{24} to infer that: if $N_1$ and $N_2$ are normal operators on ${\mathcal H}$ and $\sigma(N_1) = \sigma(N_2) = \sigma_e(N_1) = \sigma_e(N_2)$ then there is a unitary operator $U$ so that $UN_1U^* - N_2 \in {\mathcal C}_2$ and $|UN_1U^* - N_2|_2 < \varepsilon$ for a given $\varepsilon > 0$.

Also, if $T \in {\mathcal A}{\mathcal N}({\mathcal H})$ and $N$ is a normal operator with $\sigma(N) - \sigma_e(T)$ then there is a unitary operator $U: {\mathcal H} \to {\mathcal H} \oplus {\mathcal H}$ so that $(T \oplus N)U - UT$ is Hilbert--Schmidt and $|(T \oplus N)U - UT|_2 < \varepsilon$ for a given $\varepsilon > 0$,

\section{The Banach Algebras $E\Lambda(\Omega)$}
\label{sec3}

\noindent
{\bf 3.1.} We shall define here the algebras $E\Lambda(\Omega)$ and give a few of their basic properties.

If $\Omega \subset {\mathbb C}$ is a Borel set and $f \in L^{\infty}({\mathbb C},\lambda)$, with $\lambda$ denoting Lebesgue measure, let $M_f$ be the multiplication operator by $f$ on $L^2(\Omega,\lambda)$ and $Df$ be the difference quotient
\[
Df(s,t) = \frac {f(s)-f(t)}{s-t} (s \ne t)
\]
which is the class up to null-sets of a Lebesgue-measurable function on $\Omega \times \Omega$. {\em Let further}
\[
\Lambda(\Omega) = \{f \in L^{\infty}(\Omega,\lambda) \mid Df \in L^{\infty}(\Omega \times \Omega,\lambda \otimes \lambda)
\]
{\em be the subalgebra of essentially Lipschitz functions. If $T \in {\mathcal B}(L^2(\Omega,\lambda))$ let $L(T)$ be given by}
\[
L(T) = \sup\{|[M_f,T]|_2 \mid f \in \Lambda(\Omega), \|Df\|_{\infty} \le 1\}.
\]
{\em We define $E\Lambda(\Omega)$ to e the subalgebra of ${\mathcal B}(L^2(\Omega))$}
\[
E\Lambda(\Omega) = \{T \in {\mathcal B}(L^2(\Omega,\lambda)) \mid L(T) < \infty\}.
\]
{\em It is easily seen that $E\Lambda(\Omega)$ is a $*$-subalgebra of ${\mathcal B}(L^2(\Omega,\Omega))$. Even more, $E\Lambda(\Omega)$ is an involutive Banach algebra with respect to the norm $\||T\|| = \|T\| + L(T)$ and the involution is isometric $\||T\|| = \||T^*\||$.}  The proof is along standard lines and will be left to the reader.

\bigskip
\noindent
{\bf 3.2.} If $\Omega$ is specified and $w \in {\mathbb C}$, let $(e(w))(z) = \exp(i \mathop{\mathrm{Re}}(z\overline{w}))$ and let $U(w) = M_{e(w)}$, which is a unitary operator on $L^2(\Omega,\lambda)$. Also, if $\Omega$ is bounded, the multiplication operators by the functions which at $x + iy$ equal $x + iy$, $x,y$ will be denoted by $Z,X,Y$.

\bigskip
\noindent
{\bf 3.3. Proposition.} {\em If $T \in {\mathcal B}(L^2(\Omega,\lambda))$ and
\[
L_1(T) = \sup\{|w|^{-1}|[T,U(w)]|_2 \mid w \in {\mathbb C}\backslash\{0\}\}
\]
then we have $L_1(T) \le L(T) \le 2L_1(T)$ and $\||T\||_1 = \|T\| + L_1(T)$ is an equivalent Banach algebra norm on $E\Lambda(\Omega)$.

If $\Omega$ is bounded then we have
\[
L(T) = |[T,Z]|_2.
\]
}

\bigskip
\noindent
{\bf {\em Proof.}} We first establish the assertions of the proposition in case $T \in {\mathcal C}_2$. Then $T$ is given by a kernel $K \in L^2(\Omega \times \Omega,\lambda \otimes \lambda)$ and the kernel of $[M_f,T]$ is $(f(s)-f(t))K(s,t)$. The supremum of ${\mathcal C}_2$-norms of $[M_f,T]$ over all $f$ with $\|Df\|_{\infty} \le 1$ will then equal the $L^2$-norm of $(s-t)K(s,t)$, which for bounded $\Omega$ is the kernel of $[Z,T]$. On the other hand, if $f = e(w)|w|^{-1}$ we have $\|Df\|_{\infty} \le 1$, so that $L_1(T) \le L(T)$. Further, taking $w = \varepsilon w_0$, for some $w_0$ with $|w_0| = 1$ and letting $\varepsilon \downarrow 0$, the supremum of $L^2$-norms of the corresponding $(f(s)-f(t))K(s,t)$ will be the $L^2$-norm of $\mathop{\mathrm{Re}}((s-t)\overline{w}_0)K(s,t)$. The bound $L(T) \le 2L_1(T)$ is then obtained taking for instance $w_0 = 1$ and $w_0 = i$.

To deal with general $T$, we first take up the assertion that $L(T) = |[Z,T]|_2$ when $\Omega$ is bounded. Clearly it suffices to show that $L(T) \le |[Z,T]|_2$ the opposite inequality being obvious. By our results in \cite{24}, since $Z$ is a normal operator, there are finite rank projections $P_n \uparrow I$ so that $|[P_n,Z]|_2 \to 0$ as $n \to \infty$. Then if $f$ is such that $\|Df\|_{\infty} \le 1$, using the result for the Hilbert--Schmidt case, we have
\[
\begin{aligned}
|[M_f,T]|_2 &\le \limsup_{n \to \infty} |[M_f,P_nTP_n]|_2 \\
&\le \limsup_{n \to \infty} |[Z,P_nTP_n]|_2 \\
&\le \limsup_{n \to \infty} (2|[Z,P_n]|_2 \|T\| + |[Z,T]|_2) \\
&= |[Z,T]|_2.
\end{aligned}
\]

To prove the assertion about $L_1(T)$ for unbounded $\Omega$ and general $T$, we proceed along similar lines, after showing that there exist finite rank projections $P_n \uparrow 1$ so that
\[
\lim_{n \to \infty} \left( \sup_{w \in {\mathbb C}\backslash\{0\}} |[w^{-1}U(w),P_n]|_2\right) = 0.
\]
Let $\Omega_m = \{z \in \Omega \mid m - 1 \le |z| < m\}$ so that $\Omega$ is the disjoint union of the $\Omega_m$, $m \in {\mathbb N}$. On $L^2(\Omega_m,\lambda)$ we can find, by our result from \cite{24}, finite rank projections $P_{km}$ so that $P_{km} \uparrow I$ as $k \to \infty$ and $|[P_{km},Z]|_2 \le (k^2m)^{-1}$. Observe that by the result about $|[Z,T]|_2$ we proved, this gives $L(P_{km}) \le (k^2m)^{-1}$. We then define the projection $P_m$ acting on $L^2(\Omega,\lambda) = L^2(\Omega_1,\lambda) \oplus L^2(\Omega_2,\lambda) \oplus \dots$ to be $P_{m1} \oplus P_{m2} \oplus \dots \oplus P_{mm} \oplus 0 \oplus 0 \oplus \dots$ so that $P_m \uparrow I$ and $L(P_m) \le L(P_{m1}) + \dots + L(P_{mm}) \le Cm^{-1}$. Since $\|Dw^{-1}e(w)\|_{\infty} \le 1$ we have $|[w^{-1}U(w),P_m]|_2 \le Cm^{-1}$ which clearly converges to zero as $m \to \infty$ uniformly for $w \in {\mathbb C}\backslash\{0\}$. We then have for $f \in \Lambda(\Omega)$ with $\|Df\|_{\infty} \le 1$ and $T \in {\mathcal B}(L^2(\Omega,\lambda))$
\[
\begin{aligned}
|[M_f,T]|_2 &\le \limsup_{n \to \infty} |[M_f,P_nTP_n]|_2 \\
&\le \limsup_{n \to \infty} L_1(P_nTP_n) \\
&\le \limsup_{n \to \infty} (2L_1(P_n)\|T\| + L_1(T)) \\
&= L_1(T).
\end{aligned}
\]
\qed

\bigskip
\noindent
{\bf 3.4.} If $\Omega = {\mathbb C}$ the proposition provides a characterization of the algebra $E\Lambda(\Omega)$ which translates well after Fourier transform. Let ${\mathcal F}: L^2({\mathbb C},\lambda) \to L^2({\mathbb C},\lambda)$ be the unitary Fourier transform
\[
({\mathcal F} f)(w) = c \int_{{\mathbb C}} f(z)(e(-w))(z)d\lambda(z).
\]
Then ${\mathcal F} U(w_0) = V(w_0){\mathcal F}$ where $(V(w_0)g)(w) = g(w-w_0)$ and we have the following corollary.

\bigskip
\noindent
{\bf 3.5. Corollary.} {\em If $S,T \in {\mathcal B}(L^2({\mathbb C},\lambda))$ and $M(S) = \sup\{|w_0|^{-1}|S-V(w_0)SV(w_0)^*|_2 \mid w_0 \in {\mathbb C}\backslash\{0\}\}$ then we have $M({\mathcal F} T{\mathcal F}^{-1}) = L_1(T)$ and ${\mathcal F} E\Lambda({\mathbb C}){\mathcal F}^{-1} = \{S \in {\mathcal B}(L^2({\mathbb C},\lambda))| M(S) < \infty\}$.}

\bigskip
\noindent
{\bf 3.6.} If $\Omega_1 \subset \Omega_2$ let
\[
i(\Omega_2,\Omega_1): {\mathcal B}(L^2(\Omega_1,\lambda)) \to {\mathcal B}(L^2(\Omega_2,\lambda))
\]
be the inclusion homomorphism defined by $i(\Omega_2,\Omega_1)(T) = T \oplus 0$ with respect to the decomposition $L^2(\Omega_2,\lambda) = L^2(\Omega_1,\lambda) \oplus L^2(\Omega_2\backslash\Omega_1,\lambda)$. There is also a conditional expectation $\varepsilon(\Omega_1,\Omega_2): {\mathcal B}(L^2(\Omega_2,\lambda)) \to {\mathcal B}(L^2(\Omega_1,\lambda))$, $\varepsilon(\Omega_1,\Omega_2)(S) = M_{{\mathcal X}_{\Omega_1}}SM_{{\mathcal X}_{\Omega_1}} \mid L^2(\Omega_1,\lambda)$ where ${\mathcal X}_{\Omega_1}$ is the indicator function of the subset $\Omega_1$ of $\Omega_2$. It is easily checked that the Banach algebras $E\Lambda(\Omega)$ behave well with respect to the $i(\Omega_2,\Omega_1)$ and $\varepsilon(\Omega_2,\Omega_1)$.

\bigskip
\noindent
{\bf 3.7. Proposition.} {\em If $\Omega_1 \subset \Omega_2$ then we have
\[
i(\Omega_2,\Omega_1)(E\Lambda(\Omega_1)) \subset E\Lambda(\Omega_2)
\]
and the inclusion is isometric with respect to the $\||\cdot\||$-norms and also with respect to the $\||\cdot\||_1$-norms and $L(\cdot)$ and $L_1(\cdot)$ are preserved. We also have $\varepsilon(\Omega_1,\Omega_2)(E\Lambda(\Omega_2)) = E\Lambda(\Omega_1)$ and $\varepsilon(\Omega_1,\Omega_2)$ is contractive both in the $\||\cdot\||$-norms and in the $\||\cdot\||_1$-norms and we have $\varepsilon(\Omega_1,\Omega_2)i(\Omega_2,\Omega_1)(T) = T$.}

\bigskip
\noindent
{\bf 3.8.} {\em We define the Banach subalgebra $E\Lambda(\Omega)_0 \subset E\Lambda(\Omega)$ to be the closure in $E\Lambda(\Omega)$ of $\bigcup \{i(\Omega,\Omega_1)E\Lambda(\Omega_1) \mid \Omega_1 \subset \Omega,\Omega,$ bounded Borel set\}. Equivalently $E\Lambda(\Omega)_0$ is the closure in $E\Lambda(\Omega)$ of $\bigcup_{r > 0} i(\Omega,\Omega \cap r{\mathbb D})E\Lambda(\Omega \cap r{\mathbb D})$ where ${\mathbb D}$ is the unit disk.}

\bigskip
\noindent
{\bf 3.9. Proposition.} {\em $E\Lambda(\Omega)_0$ is an ideal in $E\Lambda(\Omega)$. If ${\mathcal X}_{\Omega \cap n{\mathbb D}}$ is the indicator function of $n{\mathbb D} \cap \Omega$ as a subset of $\Omega$ and $M_n = M_{{\mathcal X}_{\Omega \cap n{\mathbb D}'}}$ then $(M_n)_{n \ge 1}$ is an approximate unit of $E\Lambda(\Omega)_0$.}

\bigskip
\noindent
{\bf {\em Proof.}} Since $\||M_n\|| = \|M_n\| = 1$ and $M_nx = xM_n = x$ for any $x \in \bigcup\{i(\Omega,\Omega_1)E\Lambda(\Omega_1) \mid \Omega_1 \subset \Omega,\Omega_1$ bounded Borel\} as soon as $n$ is large enough, we clearly have that $(M_n)_{n \ge 1}$ is an approximate unit of $E\Lambda(\Omega)_0$. To prove that $E\Lambda(\Omega)_0$ is a two-sided ideal in $E\Lambda(\Omega)$ it will suffice now to show that $TM_n \in E\Lambda(\Omega)_0$ and $M_nT \in E\Lambda(\Omega)_0$. Actually since we deal with involutive algebras it will suffice to show that $TM_n \in E\Lambda(\Omega)_0$ and this in turn reduces to checking that $\||(I-M_m)TM_n\|| \to 0$ as $m \to +\infty$. It is easily seen that $L(T) < \infty$ implies $(I-M_{n+1})TM_n \in {\mathcal C}_2$ and hence $\|(I-M_m)TM_n\| \le |(I-M_m)(I-M_{n+1})TM_n|_2 \to 0$ as $m \to +\infty$. Also if $K(z_1,z_2)$ is the kernel of $(I-M_{n+1})TM_n$ then $L((I-M_{n+1})TM_n) < \infty$ means $(z_1-z_2)K(z_1,z_2)$ is in $L^2(\Omega \times \Omega,\lambda \otimes \lambda)$. Then if $m > n+1$, $L((I-M_m)TM_n)$ is the $L^2$-norm of the kernel
\[
(1 - {\mathcal X}_{\Omega \cap m{\mathbb D}}(z_1))(z_1-z_2)K(z_1,z_2)
\]
which converges to zero as $m \to +\infty$. \qed

\bigskip
\noindent
{\bf 3.10.} If $\Omega$ is bounded ${\mathcal C}_2(L^2(\Omega,\lambda)) \subset E\Lambda(\Omega)$ and $\||X\|| \le (1+d)|X|_2$ where $d$ is the diameter of $\Omega$ when $X \in {\mathcal C}_2(L^2(\Omega,\lambda))$. If $\Omega$ is unbounded the ${\mathcal C}_2\Lambda(\Omega) = {\mathcal C}_2(L^2(\Omega,\lambda)) \cap E\Lambda(\Omega)$ is only a subset of ${\mathcal C}_2(L^2(\Omega,\lambda))$. Similarly $R\Lambda(\Omega)$ will denote ${\mathcal R}(L^2(\Omega,\lambda)) \cap E\Lambda(\Omega)$ where ${\mathcal R}({\mathcal H})$ stands for the finite rank operator on ${\mathcal H}$. Remark also that if $L^2\Lambda(\Omega)$ denotes functions $f \in L^2(\Omega,\lambda)$ so that $f(z)(1+|z|) \in L^2$ then the linear span of $\langle\cdot,f\rangle g$ is in ${\mathcal R}\Lambda(\Omega)$ when $f,g \in L^2\Lambda(\Omega)$. Note also that if $f \in L^{\infty}(\Omega,\lambda)$ then $\||M_f\|| = \||Mf\||_1 = \|f\|_{\infty} = \|M_f\|$ since $L(M_f) = 0$ and $ML^{\infty}(\Omega) = \{M_f \mid f \in L^{\infty}(\Omega,\lambda)\} \subset E\Lambda(\Omega)$.

The following lemma records a consequence of the diagonalizability $\mathop{\mathrm{mod}} {\mathcal C}_2$ of normal operators, which appeared in the last part of the proof of Proposition~3.3.

\bigskip
\noindent
{\bf 3.11. Lemma.} {\em In $E\Lambda(\Omega)$ there are finite rank projections $P_n$, so that $P_n \uparrow I$ and
\[
\lim_{n \to \infty} L(P_n) = 0.
\]
Moreover we have $P_n \in i(\Omega,\Omega \cap n{\mathbb D})E\Lambda(\Omega \cap n{\mathbb D})$ and $[P_n,M_{{\mathcal X}_{\Omega \cap m{\mathbb D}}}] = 0$ for all $m \in {\mathbb N}$.}

\bigskip
We will also find it useful to have the following technical lemma when $\Omega$ is unbounded.

\bigskip
\noindent
{\bf 3.12. Lemma.} {\em Let $M_n = M_{{\mathcal X}_n} \in ML^{\infty}(\Omega,\lambda)$ where ${\mathcal X}_n$ is the indicator function of $\Omega \cap n{\mathbb D}$ as a subset of $\Omega$ and let $T \in E\Lambda(\Omega)$. Then we have $L(T-M_nTM_n) \to 0$ as $n \to \infty$.}

\bigskip
\noindent
{\bf {\em Proof.}} If $\Omega_n = \Omega \cap n{\mathbb D}$, then we have $M_nTM_n = i(\Omega,\Omega_n) \varepsilon(\Omega_n,\Omega)(T)$. With $T_n$ denoting $\varepsilon(\Omega_n,\Omega)(T)$ and $X_n$ denoting $i(\Omega,\Omega_n)([Z,T_n])$ we have the following martingale properties. If $m \ge n$ then $M_nX_mM_n = X_n$ and $|X_n|_2 = L(T_n) \le L(T)$. Hence, if $X$ is a weak limit of some subsequence of the $X_m$'s as $m \to \infty$ we will have $|X|_2 < \infty$ and $X_n = M_nXM_n$. Thus if $m \ge n$
\[
\begin{aligned}
L(M_mTM_m - M_nTM_n) &= L(\varepsilon(\Omega_m,\Omega)(M_mTM_m-M_nTM_n)) \\
&= |[Z,\varepsilon(\Omega_m,\Omega)(M_mTM_m-M_nTM_n)]|_2 \\
&= |X_m - X_n|_2.
\end{aligned}
\]
Since $M_mTM_m$ converges weakly to $T$ and $X_m$ converges in $2$-norm to $X$ as $m \to \infty$, we infer
\[
L(T-M_nTM_n) \le \sup_{m \ge n} L(M_mTM_m-M_nTM_n) = \sup_{m \ge n} |X_m-X_n|_2 = |X_m-X_n|_2.
\]
The assertion of the lemma follows from
\[
|X-X_n|_2 = |X-M_nXM_n|_2 \to 0
\]
as $n \to \infty$. \qed

\bigskip
\noindent
{\bf 3.13.} {\em We define ${\mathcal K}\Lambda(\Omega) = \{T \in E\Lambda(\Omega) \mid T$ compact\}. Clearly, ${\mathcal K}\Lambda(\Omega)$ is a closed ideal in $E\Lambda(\Omega)$.}

\bigskip
\noindent
{\bf 3.14. Proposition.} {\em The ideal ${\mathcal K}\Lambda(\Omega)$ of $E\Lambda(\Omega)$ has an approximate unit $(P_n)_{n \ge 1}$ where $P_n$'s are self-adjoint projections with the properties outlined in Lemma~$3.11$. In particular $\bigcup_{n \ge 1} P_n{\mathcal B}(L^2(\Omega,\lambda))P_n$ is a dense subalgebra in ${\mathcal K}\Lambda(\Omega)$ in $\||\cdot\||$-norm.}

\bigskip
\noindent
{\bf {\em Proof.}} If $T \in {\mathcal K}\Lambda(\Omega)$ then with the notation in Lemma~3.12 we actually have $\||T-M_nTM_n\|| \to 0$ as $n \to \infty$ in view of the lemma and of the compactness of $T$ which gives $\|T-M_nTM_n\| \to 0$. In view of the involution, the proof reduces to showing that $\||T-P_mT\|| \to 0$ as $m \to \infty$ where $P_m$ are the projections in Lemma~3.11 and $T \in {\mathcal K}\Lambda(\Omega)$ satisfies $T = M_nTM_n$ for some fixed $n$.

Clearly $T$ being compact we have $\|T-P_mT\| \to 0$ as $m \to \infty$.

On the other hand if $m \ge n$, $T-P_mT = i(\Omega,\Omega \cap n{\mathbb D})(T'-P'_mT')$ where $T' = \varepsilon(\Omega \cap n{\mathbb D},\Omega)(T)$ satisfies $i(\Omega,\Omega \cap n{\mathbb D})(T') = T$ and $P'_m = \varepsilon(\Omega \cap n{\mathbb D},\Omega)(P_m) = \varepsilon(\Omega \cap n{\mathbb D},\Omega)(P_mM_n)$ is a projection. We have
\[
\begin{aligned}
L(T-P_mT) &= L(T'-P'_mT') \\
&= |[Z,(I-P'_m)T']|_2 \\
&\le L(I-P'_m)\|T\| + |(I-P'_m)[Z,T']|_2 \to 0
\end{aligned}
\]
since $L(P'_m) \le L(P_m) \to 0$ and $[Z,T'] \in {\mathcal C}_2$, $P'_m \uparrow I$.

The remaining assertion follows from the fact that $P_n$ is an approximate unit once we remark that $P_n{\mathcal B}(L^2(\Omega,\lambda))P_n = P_nE\Lambda(\Omega)P_n = P_n{\mathcal K}\Lambda(\Omega)P_n$ because $P_n = M_nP_nM_n$.\qed

\bigskip
\noindent
{\bf 3.15. Proposition.} {\em The unit ball of $E\Lambda(\Omega)$ in $\||\cdot\||$-norm or $\||\cdot\||_1$-norm is closed in the weak operator topology and hence is weakly compact. Moreover, $E\Lambda(\Omega)$ is inverse-closed as a subalgebra of ${\mathcal B}(L^2(\Omega))$ and also closed under $C^{\infty}$-functional calculus for normal elements. In particular if $T \in E\Lambda(\Omega)$ has bounded inverse and $T = V|T|$ is its polar decomposition, then $V,|T|$ are in $E\Lambda(\Omega)$.}

\bigskip
The proof is an exercise along standard lines and will be omitted.

\bigskip
\noindent
{\bf 3.16.} {\em We shall denote by $(E/{\mathcal K})\Lambda(\Omega)$ the quotient-Banach algebra $E\Lambda(\Omega)/{\mathcal K}\Lambda(\Omega)$ and by $p: E\Lambda(\Omega) \to (E/{\mathcal K})\Lambda(\Omega)$ the canonical surjection.

Remark also that we have ${\mathcal K}\Lambda(\Omega) \subset E\Lambda(\Omega)_0$ since the dense subalgebra of ${\mathcal K}\Lambda(\Omega)$ appearing in Proposition~$3.14$ is in $E\Lambda(\Omega)_0$.
The quotient $E\Lambda(\Omega)_0/{\mathcal K}\Lambda(\Omega)$ will also be denoted $(E_0/{\mathcal K})\Lambda(\Omega)$.}

\bigskip
\noindent
{\bf 3.17. Proposition.} {\em Given $n \in {\mathbb N}$ there are $U_k \in E\Lambda(\Omega)$, $1 \le k \le n$, such that $U: L^2(\Omega,\lambda) \to L^2(\Omega,\lambda) \otimes {\mathbb C}^n$ defined by $Uh = \sum_k U_kh \otimes e_k$ is a unitary operator. In particular we have $UE\Lambda(\Omega)U^* = E\Lambda(\Omega) \otimes M_n$ and $T \to UTU^*$ is a spatial isomorphism of $E\Lambda(\Omega)$ and $E\Lambda(\Omega) \otimes {\mathcal M}_n$. Additionally we also have that
\[
UE\Lambda({\mathbb C})_0U^* = E\Lambda({\mathbb C})_0 \otimes {\mathcal M}_n.
\]
}

\bigskip
\noindent
{\bf {\em Proof.}} The existence of $U$ is a consequence of our results on normal operator $\mathop{\mathrm{mod}} {\mathcal C}_2$ (\cite{24}, $2.8$). There will be some additional technicalities due to the fact that $\Omega$ may be unbounded. From (\cite{24}, $2.8$) we get the existence of unitary operators $V_m: L^2(\Omega_m,\lambda) \to L^2(\Omega_m,\lambda) \otimes {\mathbb C}^n$, so that $|V_mZ - (Z \otimes I_n)V_m|_2 < 2^{-m-1}$, where $\Omega_m = \Omega \cap ((m+1){\mathbb D}\backslash m{\mathbb D})$. If $V_mh = \sum_k V_{mk}h \otimes e_k$, we have $|[V_{mk},Z]|_2 < 2^{-m-1}$ and hence $L(U_k) < 1$ where $U_k = \bigoplus_{m \ge 0} V_{mk}$, so that if $Uh = \sum_k U_kh \otimes e_k$ we will have that $U$ is unitary and $U_k \in E\Lambda(\Omega)$, $1 \le k \le n$. It follows that $UE\Lambda(\Omega)U^* \subset E\Lambda(\Omega) \otimes {\mathcal M}_n$ and $U^*(E\Lambda(\Omega) \otimes {\mathcal M}_n)U \subset E\Lambda(\Omega)$, which implies that $UE\Lambda(\Omega)U^* = E\Lambda(\Omega) \otimes {\mathcal M}_n$ and that $T \to UTU^*$ is a spatial isomorphism of $E\Lambda(\Omega)$ and $E\Lambda(\Omega) \otimes {\mathcal M}_n$.

For the last assertion to be proved, note that the operator $U$ which we constructed, satisfies
\[
U(i(\Omega,\Omega \cap n{\mathbb D}))E\Lambda(\Omega \cap n{\mathbb D}))U^* = (i(\Omega,\Omega \cap n{\mathbb D})E\Lambda(\Omega \cap n{\mathbb D})) \otimes {\mathcal M}_n.
\]
The assertion then follows from the density of $\bigcup_{n \ge 1} i(\Omega,\Omega \cap n{\mathbb D})E\Lambda(\Omega \cap n{\mathbb D})$ in $E\Lambda(\Omega)_0$. \qed

\bigskip
\noindent
{\bf 3.18.} Along similar lines with $3.17$ one can show that $E\Lambda(\Omega)$ is a huge algebra. For instance, since $Z$ and $Z \otimes I_{{\mathcal H}}$ are unitarily equivalent $\mathop{\mathrm{mod}} {\mathcal C}_2$ and since $I \otimes {\mathcal B}({\mathcal H})$ is in the commutant of $Z \otimes I_{{\mathcal H}}$, (${\mathcal H}$ a separable Hilbert space), one infers that $E\Lambda(\Omega)$ contains a subalgebra spatially isomorphic to $I \otimes {\mathcal B}({\mathcal H})$.

In the remainder of this section we exhibit a few special operators which are in $E\Lambda(\Omega)$.

\bigskip
\noindent
{\bf 3.19. Proposition.} {\em Let $\Omega$ be a bounded open set and let $A^2(\Omega)$ be the Bergman space of square-integrable analytic functions. Assume moreover that the rational functions with poles in ${\mathbb C}\backslash\overline{\Omega}$ are dense in $A^2(\Omega)$. Then we have $P_{\Omega} \in E\Lambda(\Omega)$, where $P_{\Omega}$ is the orthogonal projection of $L^2(\Omega,\lambda)$ onto the subspace $A^2(\Omega)$.}

\bigskip
\noindent
{\bf {\em Proof.}} This is a consequence of the Berger--Shaw inequality (see for instance \cite{20} p.~128 Theorem~1.3). Indeed $T = Z \mid A^2(\Omega)$ is a subnormal operator and the constant function $1$ is a rationally cyclic vector for $T$. The Berger--Shaw inequality then gives $\mathop{\mathrm{Tr}}[T^*,T] < \infty$. With the simplified notation $P = P_{\Omega}$, we have
\[
\mathop{\mathrm{Tr}}[PZ^*P,PZP] < \infty.
\]
Since $(I-P)ZP = 0$ and $[Z^*,Z] = 0$ this gives
\[
[PZ^*P,PZP] = PZ(I-P)Z^*P
\]
and hence
\[
[P,Z] = PZ(I-P) \in {\mathcal C}_2.
\]
\qed

\bigskip
\noindent
{\bf 3.20.} The Hilbert-transform singular integral operator on ${\mathbb C}$ (\cite{19},\cite{22})
\[
Hf(\zeta) = \lim_{\varepsilon \downarrow 0} \int_{|z-\zeta| > \varepsilon} \frac {f(z)}{(\zeta-z)^2} d\lambda(z)
\]
is a bounded operator on $L^2({\mathcal C},\lambda)$ and hence also its compression $H_{\Omega}$ to $L^2(\Omega,d\lambda)$, where $\Omega$ is bounded, is a bounded operator. Then also $T_{\Omega} = [Z,H_{\Omega}]$ is a bounded operator and
\[
T_{\Omega}f(z) = \lim_{\varepsilon \downarrow 0} \int_{|z-\zeta| > \varepsilon} \frac {f(z)}{\zeta-z} d\lambda(z).
\]
We have $[Z,T_{\Omega}] = \langle\cdot,1\rangle1$ where $1$ denotes the constant function equal to $1$. Since $[Z,T_{\Omega}]$ is rank one, we have $T_{\Omega} \in E\Lambda(\Omega)$. Since $z^{-1}$ is not in $L^2({\mathbb D},\lambda)$, $T_{\Omega}$ is not in ${\mathcal C}_2$. It can be shown that $T_{\Omega} \in {\mathcal C}_2^+$ (the ideal of compact operators with singular numbers $s_n = O(n^{-1/2})$). Also clearly the linear span of operators of the form $M_fT_{\Omega}M_g$ gives operators $K$ in $E\Lambda(\Omega)$ which are in ${\mathcal C}_2^+$ and the commutators of which $[Z,K]$ are dense in ${\mathcal C}_2(L^2(\Omega,\lambda))$.

\section{About the $K$-theory of $E\Lambda(\Omega)$}
\label{sec4}

\noindent
{\bf 4.1.} Passing via almost normal operators, the Pincus $g$-function gives a homomorphism of the $K_0$-group of $E\Lambda(\Omega)$ to $L^1$-functions. We shall prove that the Conjecture~3 about almost normal operators (see $2.5$) implies that this homomorphism completely determines the group $K_0(E\Lambda({\mathbb C})_0)$. Conversely, assuming Conjecture~4, we will show that such a result about the $K$-theory of $E\Lambda({\mathbb C})_0$ implies Conjecture~3.

We begin with some technical facts.

\bigskip
\noindent
{\bf 4.2. Lemma.} {\em If $F = F^2 \in E\Lambda(\Omega)$ and $P = P^* = P^2 \in {\mathcal B}(L^2(\Omega,\lambda))$ is the orthogonal projection onto $F(L^2(\Omega,\lambda))$ then $P \in E\Lambda(\Omega)$ and $P$ and $F$ have the same class in $K_0$.}

\bigskip
\noindent
{\bf {\em Proof.}} The orthogonal projection $P$ is equal to $\psi(FF^*)$ for some $C^{\infty}$-function $\psi$. Hence $P \in E\Lambda(\Omega)$ is a consequence of Proposition~3.15 and $tP + (1-t)F$, $t \in [0,1]$ is a continuous path of projections, so $[P]_0 = [F]_0$. \qed

\bigskip
\noindent
{\bf 4.3. Lemma.} {\em Let $P \in E\Lambda(\Omega)$ be a self-adjoint projection, which is not finite rank and assume $\Omega$ is bounded. Then we have
\[
PZP \in {\mathcal A}{\mathcal N}(L^2(\Omega,\lambda).
\]
}

\bigskip
\noindent
{\bf {\em Proof.}} We have
\[
[PZ^*P,PZP] = PZ(I-P)Z^*P - PZ^*(I-P)ZP \in {\mathcal C}_1
\]
since $(I-P)ZP = (I-P)[Z,P] \in {\mathcal C}_2$ and $PZ(I-P) = [P,Z](I-P) \in {\mathcal C}_2$. \qed

\bigskip
\noindent
{\bf 4.4. Proposition.} {\em Assume $\Omega$ is bounded. For every $\alpha \in K_0(E\Lambda(\Omega))$ there is a self-adjoint projection $P \in E\Lambda(\Omega)$, not of finite rank, so that $[P]_0 = \alpha$. The Pincus $g$-function $g_{PZP}$ depends only on $\alpha$ (i.e., not on the choice of $P$). Moreover, the map $K_0 \to L^1({\mathbb C},\lambda)$ which associates to a class $\alpha$ the $L^1$-function $g_{PZP}$ is a homomorphism.}

\bigskip
\noindent
{\bf {\em Proof.}} The existence of a unitary ``Cuntz $n$-tuple'' $U_1,\dots,U_n$ in $E\Lambda(\Omega)$, which was shown in Proposition~3.17, implies that for a projection $Q \in {\mathcal M}_n(E\Lambda(\Omega))$ there is a projection $P \in E\Lambda(\Omega)$ with $[P]_0 = [Q]_0$ and that $[I]_0 = 0$, so that $-[Q]_0 = [I-P]_0$. Hence $K_0(E\Lambda(\Omega))$ consists of classes of idempotents in $E\Lambda(\Omega)$ and these can be chosen to be self-adjoint by Lemma~4.2.

Again using Proposition~3.15 and Proposition~3.17 the fact that the map $\alpha \to g_{PZP}$ is a well-defined homomorphism is a consequence of the following two facts: a) if $P \in E\Lambda(\Omega)$ is a self-adjoint projection and $W \in E\Lambda(\Omega)$ is unitary, then $g_{(WPW^*)Z(WPW^*)} = g_{PZP}$ and b) if $P_1,P_2 \in E\Lambda(\Omega)$ are self-adjoint projections and $P_1P_2 = 0$, then $g_{P_1ZP_1} + g_{P_2ZP_2} = g_{(P_1+P_2)Z(P_1+P_2)}$.

To show that a) holds, remark that $g_{WPW^*ZWPW^*} = g_{PW*ZWP}$ by unitary equivalence and $PW^*ZWP-PZP \in {\mathcal C}_2$. Moreover, in view of the argument in $2.7$ we have $k_2(PZP) = 0$, $k_2(PW^*ZWP) = 0$ and we can then use $2.4$ to get that $g_{PZP} = g_{PW^*ZWP}$.

Assertion b) is proved by the same kind of combination of facts. By the argument of $2.7$, we have
\[
k_2(P_1ZP_1) = k_2(P_2ZP_2) = k_2((P_1+P_2)Z(P_1+P_2)) = 0.
\]
We then remark that
\[
P_1ZP_1 + P_2ZP_2 - (P_1+P_2)Z(P_1+P_2) \in {\mathcal C}_2
\]
and we can then use $2.4$ to get
\[
g_{(P_1+P_2)Z(P_1+P_2)} = g_{P_1ZP_1 + P_2ZP_2} = g_{P_1ZP_1} + g_{P_2ZP_2},
\]
where we used the fact that
\[
k_2(P_1ZP_1 + P_2ZP_2) = k_2(P_1ZP_1 \oplus P_2ZP_2) = 0.
\]
\qed

\bigskip
\noindent
{\bf 4.5.} {\em The homomorphism $K_0(E\Lambda(\Omega)) \to L_{rc}^1({\mathbb C},\lambda)$, constructed in Proposition~$4.4$, will be denoted by $\Gamma(\Omega)$ or simply $\Gamma$, when the bounded set $\Omega$ is not in doubt $(L_{rc}^1({\mathbb C},\lambda)$ being the $L^1$-space of real-valued functions with compact support).

We shall also denote by ${\mathcal A}{\mathcal N}{\mathcal D}({\mathcal H})$ the almost normal operators for which Conjecture~$4$ (see $2.5$) holds. We shall call such almost-normal operators dilatable. It is easily seen that this is equivalent to the fact that the almost-normal operator is a Hilbert--Schmidt perturbation of an almost-normal operator which is a compression $PNP$ of a normal operator $N$ by a projection $P$ so that $[P,N] \in {\mathcal C}_2$.

In $2.7$ we showed that if $T \in {\mathcal A}{\mathcal N}{\mathcal D}({\mathcal H})$ then $k_2(T) = 0$.}

\bigskip
Next we will give a few simple facts about $K$-theory for some of the algebras related to $E\Lambda(\Omega)$ and get some variants of the homomorphism $\Gamma$.

\bigskip
\noindent
{\bf 4.6.} If $\Omega_1 \subset \Omega_2$ are bounded Borel sets, then it is immediate from the construction of $\Gamma$ that
\[
\Gamma(\Omega_2) \circ (i(\Omega_2,\Omega_1))_* = \Gamma(\Omega_1).
\]
In view of 3.8, $E\Lambda({\mathbb C})_0$ is the inductive limit of the $E\Lambda(\Omega)$ with $\Omega$ bounded (the inclusion will be denoted $i_0({\mathbb C},\Omega)$). {\em Then $K_0(E\Lambda({\mathbb C})_0)$ is the inductive limit of the $K_0(E\Lambda(\Omega))$, with bounded $\Omega$, and there is a homomorphism
\[
\Gamma_{\infty}: K_0(E\Lambda({\mathbb C})_0) \to L_{rc}^1({\mathbb C},\lambda)
\]
so that
\[
\Gamma_{\infty} \circ (i_0({\mathbb C},\Omega))_* = \Gamma(\Omega).
\]
}

\bigskip
\noindent
{\bf 4.7. Lemma.} {\em We have $K_0({\mathcal K}\Lambda(\Omega)) \cong {\mathbb Z}$, $K_1({\mathcal K}\Lambda(\Omega)) = 0$, for any $\Omega$ (not of measure $0$), the isomorphism for $K_0$ being given by the trace on $B(L^2(\Omega,
\lambda))$. Moreover we have isomorphisms
\[
\begin{aligned}
K_0(E\Lambda(\Omega)) &\overset{p_*}{\longrightarrow} K_0((E/{\mathcal K})\Lambda(\Omega)) \\
K_0(E\Lambda({\mathbb C})_0) &\overset{p_*}{\longrightarrow} K_0((E_0/{\mathcal K})\Lambda({\mathbb C})).
\end{aligned}
\]
}

\bigskip
\noindent
{\bf {\em Proof.}} The assertions about the $K$-theory of ${\mathcal K}\Lambda(\Omega)$ are a consequence of the last assertion in Proposition~3.14.

To get the isomorphisms between $K_0$-groups of $E\Lambda(\Omega)$ and $(E/{\mathcal K})\Lambda(\Omega)$ and respectively $E\Lambda({\mathbb C})_0$ and $(E_0/{\mathcal K})\Lambda({\mathbb C})$ we use the $6$-term $K$-theory exact sequences associated with
\[
0 \to {\mathcal K}\Lambda(\Omega) \to E\Lambda(\Omega) \to (E/{\mathcal K})\Lambda(\Omega) \to 0
\]
and
\[
0 \to {\mathcal K}\Lambda({\mathbb C}) \to E\Lambda({\mathbb C})_0 \to (E_0/{\mathcal K})\Lambda({\mathbb C}) \to 0.
\]
Since $K_1({\mathcal K}\Lambda(\Omega)) = 0$ we have that the homomorphisms $p_*$ are surjective. The injectivity of the $p_*$ means to show the connecting homomorphisms $K^1 \to K^0$ are surjective. This is easily seen to be the case if we can prove $E\Lambda(\Omega)$ and $E\Lambda({\mathbb C})_0 + {\mathbb C} I$ contain a Fredholm operator of index $1$. If $\Omega$ is a Fredholm operator of index $1$, $T \in {\mathcal B}(L^2(\Omega,\lambda))$ so that $[T,Z] \in {\mathcal C}_2$. This in turn follows from the easily seen fact that $Z$ is unitarily equivalent to $Z \oplus \mu I_{{\mathcal H}} + K$, where ${\mathcal H}$ is some infinite-dimensional Hilbert space, $\mu \in \sigma(Z)$ and $K \in {\mathcal C}_2$. For $E\Lambda({\mathbb C})_0 + {\mathbb C} I$ we can use the Fredholm operator $T \in E\Lambda(\Omega)$ and consider $T \oplus I_{L^2({\mathbb C}\backslash\Omega,\lambda)} \in E\Lambda({\mathbb C})_0 + {\mathbb C} I$. \qed

\bigskip
\noindent
{\bf 4.8.} {\em In view of Lemma~$4.7$ we infer for bounded $\Omega$ the existence of homomorphisms
\[
{\tilde \Gamma}(\Omega): K_0((E/{\mathcal K})\Lambda(\Omega)) \to L_{rc}^1(\Omega,\lambda)
\]
and
\[
{\tilde \Gamma}_{\infty}: K_0((E_0/{\mathcal K})\Lambda({\mathbb C})) \to L_{rc}^1({\mathbb C},\lambda)
\]
so that
\[
\begin{aligned}
{\tilde \Gamma}(\Omega) \circ p_* &= \Gamma(\Omega) \text{ and} \\
{\tilde \Gamma}_{\infty} \circ p_* &= \Gamma_{\infty}.
\end{aligned}
\]
}

\bigskip
\noindent
{\bf 4.9. Fact.} {\em The following assertions are equivalent.}
\begin{itemize}
\item[(i)] {\em Conjecture~$3$ is true.}
\item[(ii)] {\em Conjecture~$4$ is true and $\Gamma_{\infty}$ is an isomorphism.}
\item[(iii)] {\em Conjecture~$4$ is true and $\Gamma_{\infty}$ is injective.}
\end{itemize}

\bigskip
\noindent
{\bf {\em Proof.}} Since (ii)~$\Rightarrow$~(iii) it will be sufficient to show that (i)~$\Rightarrow$~(ii) and (iii)~$\Rightarrow$~(i).

(i)~$\Rightarrow$~(ii). Remark first that Conjecture~3 implies Conjecture~4. Indeed, if $T \in {\mathcal A}{\mathcal N}({\mathcal H})$ we can find $S_1 \in {\mathcal A}{\mathcal N}({\mathcal H})$ so that $g_{S_1} = -g_T$ (see \cite{20} for instance). Then Conjecture~3 implies that there is a normal operator $N_1$ so that $T \oplus S_1 \oplus N_1 - N \in {\mathcal C}_2$ where $N$ is a normal operator. Thus we can take $S = S_1 \oplus N_1$ and then $S \in {\mathcal A}{\mathcal N}$ and $T \oplus S$ is equal $N \mathop{\mathrm{mod}} {\mathcal C}_2$, which is the assertion of Conjecture~4 for $T$.

To show $\Gamma_{\infty}$ is surjective consider $g \in L_{rc}^1({\mathbb C},\lambda)$. By the work of Carey--Pincus there is $T \in {\mathcal A}{\mathcal N}({\mathcal H}_1)$ so that $g_T = g$. By Conjecture~4 and the fact that it implies Conjecture~1 we see that $T$ can be chosen to be $QN \mid Q{\mathcal H}$ where $N$ is a normal operator and $Q$ an orthogonal projection, so that $[Q,N] \in {\mathcal C}_2$. We may also assume $\sigma(N) = n\overline{{\mathbb D}}$ for some $n \in {\mathbb N}$. Then by our results on normal operators $\mathop{\mathrm{mod}} {\mathcal C}_2$, there is a unitary operator $U: {\mathcal H} \to L^2(n{\mathbb D},\lambda)$ so that $ZU - UN \in {\mathcal C}_2$. Then taking $P = UQU^*$, we have $PZP - UQNQU^* \in {\mathcal C}_2$ and hence $g_{PZP} = g_{QNQ} = g$ so that $\Gamma(n{\mathbb D})[P]_0 = g$. Clearly then $\Gamma_{\infty}([i_0({\mathbb C},n{\mathbb D})(P)]_0) = g$.

To prove that assuming Conjecture~3 holds, $\Gamma_{\infty}$ is injective, let $\alpha \in K_0(E\Lambda({\mathbb C})_0)$ be so that $\Gamma_{\infty}(\alpha) = 0$. Using 4.6 and Proposition~4.4 there is a self-adjoint projection $P \in E\Lambda(n{\mathbb D})$ for some $n \in {\mathbb N}$, so that $(i_0({\mathbb C},n{\mathbb D}))_*[P]_0 = \alpha$ and $\Gamma(n{\mathbb D})([P]_0) = \Gamma_{\infty}(\alpha) = 0$. Hence $g_{PZP} = 0$. Then Conjecture~3 gives that there is $m \ge n$ and there are normal operators $N$ and $N_1$ with $\sigma(N) = \sigma(N_1) = m\overline{{\mathbb D}}$ so that
\[
N - PZ \mid PL^2(n{\mathbb D},\lambda) \oplus N_1 \in {\mathcal C}_2.
\]
Since we will use the operators $Z$ in $E\Lambda(n{\mathbb D})$ and $E\Lambda(m{\mathbb D})$ simultaneously, we shall denote them here by $Z_n$ and $Z_m$. Clearly, we may use a unitary equivalence and a ${\mathcal C}_2$-perturbation to choose $N_1$. Similarly $N$ can be chosen unitarily equivalent to $Z_m$. Thus, we get a unitary operator
\[
U: PL^2(n{\mathbb D},\lambda) \oplus L^2(m{\mathbb D},\lambda) \to L^2(m{\mathbb D},\lambda)
\]
so that $Z_mU - U(PZ_n \mid PL^2(n{\mathbb D},\lambda) \oplus Z_m) \in {\mathcal C}_2$. This means that $U$ gives rise to a partial isometry $W \in {\mathcal B}(L^2(m{\mathbb D},\lambda) \oplus L^2(m{\mathbb D},\lambda))$ so that $W^*W = i(m{\mathbb D},n{\mathbb D})(P) \oplus I$ and $WW^* = 0 \oplus I$ with the property that $[W,Z_m \oplus Z_m] \in {\mathcal C}_2$. Then we have $W \in {\mathcal M}_2(E\Lambda(m{\mathbb D}))$. This gives $i(m{\mathbb D},n{\mathbb D})_*[P]_0 + [I]_0 = [I]_0$ in $K_0(E\Lambda(m{\mathbb D},\lambda))$, so that $[i(m{\mathbb D},n{\mathbb D})(P)]_0 = 0$. But then we must have $\alpha = [i_0({\mathbb C},n{\mathbb D})(P)]_0 = i_0({\mathbb C},m{\mathbb D})_*[i(m{\mathbb D},n{\mathbb D})(P)]_0 = 0$.

(iii)~$\Rightarrow$~(i). Assume (iii) holds and let $T_1,T_2 \in {\mathcal A}{\mathcal N}({\mathcal H})$ with $g_{T_1} = g_{T_2}$. Since Conjecture~4 is part of the assumption (iii) we have $T_1,T_2 \in {\mathcal A}{\mathcal N}{\mathcal D}({\mathcal H})$. This implies there are self-adjoint projection $P_1,P_2 \in E\Lambda(n{\mathbb D})$ for some $n \in {\mathbb N}$, so that $T_j$ is unitarily equivalent to a ${\mathcal C}_2$-perturbation of $P_jZ \mid P_jL^2(n{\mathbb D},\lambda)$, $j = 1,2$. Moreover, we have $\Gamma(n{\mathbb D})[P_1]_0 = \Gamma(n{\mathbb D})[P_2]_0$ because $g_{T_1} = g_{T_2}$. It follows that $\Gamma_{\infty}([i_0({\mathbb C},n{\mathbb D})(P_1)]_0) = \Gamma_{\infty}([i_0({\mathbb C},n{\mathbb D})(P_2)]_0)$ so that by (iii) we have $i_0({\mathbb C},n{\mathbb D})_*[P_1]_0 = i_0({\mathbb C},n{\mathbb D})[P_2]_0$. Since $E\Lambda({\mathbb C})_0$ is the inductive limit of the $E\Lambda(m{\mathbb D})$ we infer that $[i(m{\mathbb D},n{\mathbb D})(P_1)]_0 = [i(m{\mathbb D},n{\mathbb D})(P_2)]_0$ for some $m \ge n$. Hence there is a unitary equivalence in ${\mathcal M}_{p+q+1}(E\Lambda(m{\mathbb D}))$ between the $Q_j = i(m{\mathbb D},n{\mathbb D})P_j \oplus I \oplus \dots \oplus I \oplus 0 \oplus \dots \oplus 0$, $j = 1,2$ (there are $p$ summands $I$ and $q$ summands $0$). Indeed the equality of $K_0$-classes implies there is an invertible element intertwining $Q_1,Q_2$ and using Proposition~3.17 and Proposition~3.15 we can pass to the unitary in the polar decomposition of this invertible element of ${\mathcal M}_{p+q+1}(E\Lambda(n{\mathbb D}))$. This unitary will then commute with $Z \oplus \dots \oplus Z$ modulo ${\mathcal C}_2$ and hence will intertwine $\mathop{\mathrm{mod}} {\mathcal C}_2$ the compressions $Q_j(Z \oplus \dots \oplus Z)Q_j$, $j = 1,2$. These compressions are unitarily equivalent to 
\[
P_jZ \mid P_jL^2(n{\mathbb D},\lambda) \oplus N_j
\]
for some normal operators $N_j$, $j = 1,2$. Thus $T_j \oplus N_j$, being unitarily equivalent $\mathop{\mathrm{mod}} {\mathcal C}_2$ to these compressions, will also be unitarily equivalent $\mathop{\mathrm{mod}} {\mathcal C}_2$, which proves (i) under the assumption (iii). \qed

\bigskip
\noindent
{\bf 4.10.} In view of Lemma~4.7 and of 4.8 we have that Fact~4.9 also holds with $\Gamma_{\infty}$ replaced by ${\tilde \Gamma}_{\infty}$.

\section{Multipliers, Corona and Bidual of ${\mathcal K}\Lambda(\Omega)$}
\label{sec5}

\noindent
{\bf 5.1.} We shall consider bounded multipliers ${\mathcal M}({\mathcal K}\Lambda(\Omega))$, that is double centralizer pairs $(T',T'')$ of bounded linear maps ${\mathcal K}\Lambda(\Omega) \to {\mathcal K}\Lambda(\Omega)$ so that $T'(x)y = xT''(y)$.

\bigskip
\noindent
{\bf 5.2. Proposition.} {\em We have ${\mathcal M}({\mathcal K}\Lambda(\Omega)) = E\Lambda(\Omega)$, that is, if $(T',T'') \in {\mathcal M}({\mathcal K}\Lambda(\Omega))$, then there is $T \in E\Lambda(\Omega)$ so that $T'(x) = xT$ and $T''(x) = Tx$.}

\bigskip
\noindent
{\bf {\em Proof.}} Let $(P_n)_{n \ge 1}$ be the approximate unit provided by Proposition~3.14 and define $K_n = T'(P_n)P_n = P_nT''(P_n)$. Clearly, the norms $\||K_n\||$ will be bounded by some constant $C$ and if $m > n$ we have
\[
\begin{aligned}
P_nK_mP_n &= P_nT'(P_m)P_mP_n \\
&= P_nT'(P_m)P_n = P_nP_mT''(P_n) \\
&= P_nT''(P_n) = K_n.
\end{aligned}
\]
Hence if $T$ is the weak limit of the $K_n$'s we shall have $P_nTP_n = K_n$. Also $L(T) \le \sup_n(L(K_n) + 2\|T\|L(P_n)) < \infty$, so that $T \in E\Lambda(\Omega)$. Moreover, we have
\[
\begin{aligned}
T'(P_n) &= w - \lim_{m \to \infty} T'(P_n)P_m \\
&= w - \lim_{m \to \infty} P_nT''(P_m) \\
&= w - \lim_{m \to \infty} P_nP_mT''(P_m) = P_nT
\end{aligned}
\]
and similarly $T''(P_n) = TP_n$. This gives $P_nT''(x) = T'(P_n)x = P_nTx$ if $x \in {\mathcal K}\Lambda(\Omega)$ and hence 
\[
T''(x) = \lim_{n \to \infty} P_n T''(x) = \lim_{n \to \infty} P_nTx = Tx.
\]
Similarly $T'(x)P_n = xTP_m$ and $T(x) = \lim_{n \to \infty} T'(x)P_n = xT$. \qed

\bigskip
\noindent
{\bf 5.3. Proposition.} {\em The involutive Banach algebra $(E/{\mathcal K})\Lambda(\Omega)$ is a $C^*$-algebra. Actually if $x \in E\Lambda(\Omega)$ the norm of $p(x)$ in $(E/{\mathcal K})\Lambda(\Omega)$ is equal to the norm of $x+{\mathcal K}$ in the Calkin algebra ${\mathcal B}/{\mathcal K}$. In particular $(E/{\mathcal K})\Lambda(\Omega)$ is isometrically isomorphic to a $C^*$-subalgebra of ${\mathcal B}/{\mathcal K}$.}

\bigskip
\noindent
{\bf {\em Proof.}} It is easily seen that all assertions follow from the equality of the norm of $p(x)$ with the norm of $x+{\mathcal K}$ in the Calkin algebra. This in turn will follow from the fact that with $(P_n)_{n\ge 1}$ denoting the approximate unit of ${\mathcal K}\Lambda(\Omega)$ in Proposition~3.14
\[
\lim_{n \to \infty} \|(I-P_n)x(I-P_n)\|
\]
equals the Calkin norm of $x+{\mathcal K}$, if we will also show that
\[
\lim_{n \to \infty} L((I-P_n)x(I-P_n)) = 0.
\]
In case $\Omega$ is bounded we indeed have
\[
L((I-P_n)x(I-P_n)) \le \lim_{n \to \infty} (2\|x\|\,|[I-P_n,z]|_2 + |(I-P_n)[Z,x](I-P_n)|_2 = 0.
\]
In case $\Omega$ is unbounded we use Lemma~3.12 and write $x = x_0 + x_1$ where $x_0 = M_m \times M_m$ with $m$ chosen so that $L(x_1) < \varepsilon$. We have
\[
\limsup_{n \to \infty} L((I-P_n)x_1(I-P_n)) \le L(x_1) < \varepsilon
\]
and since $\varepsilon > 0$ can be chosen arbitrarily small it will suffice to show that
\[
\limsup_{n \to \infty} L((I-P_n)x_0(I-P_n)) = 0.
\]
This in turn can be seen as follows. Let $Z_k$ be the multiplication operator by $z(1 \wedge k|z|^{-1})$. Then for any $y \in E\Lambda(\Omega)$ we have
\[
L(y) = \limsup_{k \to \infty} |[Z_k,y]|_2.
\]
Moreover if $k \ge m$, $[Z_k,x_0] = [Z_m,x_0]$. Hence
\[
L((I-P_n)x_0(I-P_n)) \le 2\|x_0\|L(I-P_n) + |(I-P_n)[Z_m,x_0](I-P_n)|_2 \to 0
\]
as $n \to \infty$. \qed

\bigskip
\noindent
{\bf 5.4. Remark.} {\em The $C^*$-algebra
\[
\{p(M_f) \in (E/{\mathcal K})\Lambda(\Omega) \mid f \in C\Lambda(\Omega)\},
\]
where $C\Lambda(\Omega)$ denotes the norm closure of $\Lambda(\Omega)$ in $L^{\infty}(\Omega,\lambda)$, is in the center of $(E/{\mathcal K})\Lambda(\Omega)$.}

\bigskip
Indeed, if $f \in \Lambda(\Omega)$ then $[M_f,x] \in {\mathcal C}_2\Lambda \subset {\mathcal K}\Lambda(\Omega)$ if $x \in E\Lambda(\Omega)$ so that $p(M_f)$ is in the center of $(E/{\mathcal K})\Lambda(\Omega)$. Since $\||M_f\|| = \|M_f\| = \|f\|_{\infty}$ if $f \in L^{\infty}(\Omega)$ and the center is clearly norm-closed in $(E/{\mathcal K})\Lambda(\Omega)$, the assertion follows.

\bigskip
\noindent
{\bf 5.5.} We pass to describing the dual of ${\mathcal K}\Lambda(\Omega)$ for bounded $\Omega$. Throughout ${\mathcal C}_1$ and ${\mathcal C}_2$ will stand for ${\mathcal C}_1(L^2(\Omega,\lambda))$ and respectively ${\mathcal C}_2(L^2(\Omega,\lambda))$.

\bigskip
\noindent
{\bf 5.6. Proposition.} {\em Assuming $\Omega$ is bounded, the dual of ${\mathcal K}\Lambda(\Omega)$ can be identified isometrically with $({\mathcal C}_1 \times {\mathcal C}_2)/{\mathcal N}$ where
\[
{\mathcal N} = \{([Z,H],H) \in {\mathcal C}_1 \times {\mathcal C}_2 \mid H \in {\mathcal C}_2 \text{ with } [Z,H] \in {\mathcal C}_1\}
\]
and the duality map ${\mathcal K}\Lambda(\Omega) \times ({\mathcal C}_1 \times {\mathcal C}_2) \to {\mathbb C}$ is $(T,(x,y)) = \mathop{\mathrm{Tr}}(Tx + [Z,T]y)$.}

\bigskip
\noindent
{\bf {\em Proof.}} Since $T \to T \oplus [Z,T]$ identifies ${\mathcal K}\Lambda(\Omega)$ isometrically with a closed subspace of ${\mathcal K} \oplus {\mathcal C}_2$ endowed with the norm $\|K \oplus H\| = \|K\| + |H|_2$, the dual of which is ${\mathcal C}_1 \times {\mathcal C}_2$, the proof will boil down to showing that ${\mathcal N}$ is the annihilator of
\[
\{T \oplus [Z,T] \in {\mathcal K} \oplus {\mathcal C}_2 \mid T \in {\mathcal K}\Lambda(\Omega)\}.
\]
Since the set ${\mathcal R}$ of finite rank operators is dense in ${\mathcal K}\Lambda(\Omega)$, it will be sufficient to show that ${\mathcal N}$ is the annihilator of
\[
\{R \oplus [Z,R] \in {\mathcal K} \oplus {\mathcal C}_2 \mid R \in {\mathcal R}\}.
\]
If $R \in {\mathcal R}$ and $(x,y) \in {\mathcal N}$ we have
\[
\mathop{\mathrm{Tr}}(Rx + [Z,R]y) = \mathop{\mathrm{Tr}}(R[Z,y] + [Z,R]y) = \mathop{\mathrm{Tr}}([Z,Ry]) = 0.
\]
Conversely if $(x,y) \in {\mathcal C}_1 \times {\mathcal C}_2$ is such that
\[
\mathop{\mathrm{Tr}}(Rx + [Z,R]y) = 0 \text{ for all } R \in {\mathcal R},
\]
then
\[
\mathop{\mathrm{Tr}}(R(x-[Z,y])) = 0 \text{ for all} R \in {\mathcal R}
\]
and hence $x = [Z,y]$, that is $(x,y) \in {\mathcal N}$. \qed

\bigskip
\noindent
{\bf 5.7. Lemma.} {\em Under the same assumptions and notations like in $5.6$,
\[
\{([Z,R],R) \in {\mathcal C}_1 \times {\mathcal C}_2 \mid R \in {\mathcal R}\}
\]
is dense in ${\mathcal N}$.}

\bigskip
\noindent
{\bf {\em Proof.}} Let $(x,y) \in {\mathcal N}$, that is $y \in {\mathcal C}_2$ is such that $x = [Z,y] \in {\mathcal C}_1$. Let $(P_n)_{n \ge 1}$ be self-adjoint projections of finite rank so that $P_n \uparrow I$ and $|[P_n,Z]|_2 \to 0$. Then we have $|yP_n-y|_2 \to 0$ and also
\[
\begin{aligned}
|[Z,yP_n] - [Z,y]|_1 &= |[Z,y]P_n + y[Z,P_n] - [Z,y]|_1 \\
&\le |y|_2|[Z,P_n]|_2 + |[Z,y](I-P_n)|_1 \to 0
\end{aligned}
\]
as $n \to \infty$. \qed

\bigskip
\noindent
{\bf 5.8. Proposition.} {\em If $\Omega$ is bounded, with the same notations as in Proposition~$5.6$, the dual of $({\mathcal C}_1 \times {\mathcal C}_2)/{\mathcal N}$ identifies with $E\Lambda(\Omega)$ via the duality map
\[
(T,(x,y)) \to \mathop{\mathrm{Tr}}(Tx + [Z,T]y).
\]
In particular $E\Lambda(\Omega)$ identifies with the bidual of ${\mathcal K}\Lambda(\Omega)$.}

\bigskip
\noindent
{\bf {\em Proof.}} The dual of $({\mathcal C}_1 \times {\mathcal C}_2)/{\mathcal N}$ is the orthogonal of ${\mathcal N}$ in ${\mathcal B} \oplus {\mathcal C}_2 = ({\mathcal C}_1 \times {\mathcal C}_2)^d$ (the usual duality based on the trace). Since Lemma~5.8 provides a dense subset of ${\mathcal N}$, it suffices to show that $\{T \oplus [Z,T] \in {\mathcal B} \oplus {\mathcal C}_2 \mid T \in E\Lambda(\Omega)\}$ is the orthogonal in ${\mathcal B} \oplus {\mathcal C}_2$ of $\{([Z,R],R) \in {\mathcal C}_1 \times {\mathcal C}_2 \mid R \in {\mathcal R}\}$. Indeed, if $T \oplus H \in {\mathcal B} \oplus {\mathcal C}_2$ is such that $\mathop{\mathrm{Tr}}(T[Z,R] + HR) = 0$ for all $R \in {\mathcal R}$, then $\mathop{\mathrm{Tr}}((-[Z,T] + H)R) = 0$ for all $R \in {\mathcal R}$ and hence $H = [Z,T]$, which also implies $T \in E\Lambda(\Omega)$. Clearly, also if $T \in E\Lambda(\Omega)$ and $R \in {\mathcal R}$ we have
\[
\mathop{\mathrm{Tr}}(T[Z,R] + [Z,T]R) = \mathop{\mathrm{Tr}}([Z,TR]) = 0.
\]
\qed

\section{Concluding Remarks}
\label{sec6}

\noindent
{\bf 6.1. Isomorphisms induced by bi-Lipschitz map.} Let $\Omega_1$ and $\Omega_2$ be Borel subsets of ${\mathbb C}$ and let $F: \Omega_1 \to \Omega_2$ be a map which is Lipschitz and has an inverse which is also Lipschitz (i.e., $F$ is bi-Lipschitz). Then if $\lambda_j$ is the restriction of Lebesgue measure to $\Omega_j$, the measures $F_*\lambda_1$ are $\lambda_2$ are mutually absolutely continuous with bounded Radon--Nikodym derivatives and the same holds for $(F^{-1})_*\lambda_2$ and $\lambda_1$ (\cite{16}). This gives rise to a unitary operator
\[
U(\Omega_2,\Omega_1)L^2(\Omega_1,\lambda_1) \to L^2(\Omega_2,\lambda_1)
\]
which maps $f \in L^2(\Omega_1,\lambda_1)$ to $(f \circ F^{-1}) \cdot (dF_*\lambda_1/d\lambda_2)^{1/2}$. If $g \in L^{\infty}(\Omega_2,\lambda_2)$ then
\[
U(\Omega_2,\Omega_1)^{-1}M_gU(\Omega_2,\Omega_1) = M_{g \circ F}.
\]
The map $g \to g \circ F$ gives isomorphisms of $L^{\infty}(\Omega_2,\lambda_2)$ with $L^{\infty}(\Omega_1,\lambda_1)$ and of $\Lambda(\Omega_2)$ with $\Lambda(\Omega_1)$. Further $T \to U(\Omega_2,\Omega_1)^{-1}TU(\Omega_2,\Omega_1)$ is an isomorphism of $E\Lambda(\Omega_2)$ and $E\Lambda(\Omega_1)$. This is an isomorphism of Banach algebras with involution, which however is not isometric, since its norm depends on the Lipschitz constants of $F$ and $F^{-1}$. These isomorphisms preserve finite-rank operators and hence ${\mathcal K}\Lambda(\Omega_2)$ is mapped onto ${\mathcal K}\Lambda(\Omega)$. This in turn implies there is an induced $C^*$-algebra isomorphism of $(E/{\mathcal K})\Lambda(\Omega_2)$ with $(E/{\mathcal K})\Lambda(\Omega_1)$.

In particular the group of bi-Lipschitz homeomorphisms of a Borel set $\Omega$ has automorphic actions on $E\Lambda(\Omega)$ and $(E/{\mathcal K})\Lambda(\Omega)$.

\bigskip
\noindent
{\bf 6.2.} In view of 5.6 it is a {\em natural question to ask, what is the center of $(E/{\mathcal K})\Lambda(\Omega)$}? Note that the answer to the Calkin-algebra analogue of this question, that is the determination of the center of the commutant of a separable commutative $C^*$-subalgebra of the Calkin algebra, is a particular case of our Calkin algebra bicommutant theorem (\cite{23}).

\bigskip
\noindent
{\bf 6.3. ${\mathcal K}\Lambda(\Omega)$ as a Dirichlet algebra.} The algebras ${\mathcal K}\Lambda(\Omega)$ are examples of Dirichlet algebras in the sense of non-commutative potential theory (\cite{1}, \cite{9}, \cite{10}). The Dirichlet form can be described for instance via the construction of Dirichlet forms from derivations (Theorem~4.5 in \cite{9} or Theorem~8.3 in \cite{10}). This corresponds to working with the $C^*$-algebra of compact operators ${\mathcal K} = {\mathcal K}(L^2(\Omega,\lambda))$ and its trace Tr, which is densely defined, faithful, semifinite and lower semicontinuous. The Hilbert space ${\mathcal H} = {\mathcal C}_2 \oplus {\mathcal C}_2$, where ${\mathcal C}_2 = {\mathcal C}_2(L^2(\Omega,\lambda))$ is a ${\mathcal K}-{\mathcal K}$-bimodule and ${\mathcal J}(x \oplus y) = x^* \oplus y^*$ is an isometric antilinear involution of ${\mathcal H}$ exchanging the right and left actions of ${\mathcal K}$ on ${\mathcal H}$. Clearly ${\mathcal C}_2$ identifies with $L^2({\mathcal K},\mathop{\mathrm{Tr}})$ and there is an $L^2$-closable derivation $\partial$ of ${\mathcal K}\Lambda(\Omega) \cap {\mathcal C}_2 \to {\mathcal C}_2$. The definition in case $\Omega$ is bounded, is $\partial a = [X,a] \oplus [Y,a]$. In general the definition can be given in terms of he kernel $K(z_1,z_2)$ of an element $a \in {\mathcal K}\Lambda(\Omega) \cap {\mathcal C}_2$. Then the components of $\partial a$ have kernels $(x_1-x_2)K(z_1,z_2)$ and respectively $(y_1-y_2)K(z_1,z_2)$, which are square integrable since $a \in {\mathcal K}\Lambda(\Omega)$. Also clearly viewed the domain of definition of $\partial$ as part of $L^2({\mathcal K},\mathop{\mathrm{Tr}})$, the map $\partial$ is $L^2$-closed. Moreover $\partial$ satisfies the symmetry condition ${\mathcal J} \partial a = \partial a^*$. Then the Dirichlet form ${\mathcal E}$ which is obtained as the closure ${\mathcal E}[a] = \|\partial a\|_{{\mathcal H}}^2$ is easily seen to be precisely square of the $L^2$-norm of $(z_1-z_2)K(z_1,z_2)$ which is the same as $(L(a))^2$ defined for $a \in {\mathcal K}\Lambda(\Omega)$. The Markovian semigroup $T_t$ will then act on elements $a \in {\mathcal K}\Lambda(\Omega) \cap {\mathcal C}_2$ which have kernels $K(z_1,z_2)$ as a multiplier which produces the element with kernel $e^{-t|z_1-z_2|^2}K(z_1,z_2)$. In view of the Markovianity it is easy to see that $T_t$ extends to a semigroup of completely positive contraction on ${\mathcal K}\Lambda(\Omega)$, $E\Lambda(\Omega)$ and also on ${\mathcal K}$ and ${\mathcal B}$. Moreover $T_t$ also induces a semigroup of completely positive contractions on $(E/{\mathcal K})\Lambda(\Omega)$.

\bigskip
\noindent
{\bf 6.4. Replacing ${\mathcal C}_2$ by some other ${\mathcal C}_p$.} One may wonder about the consequences of replacing the Hilbert--Schmidt class ${\mathcal C}_2$ by some other ${\mathcal C}_p$-class in the definition of $E\Lambda(\Omega)$. This would mean to consider operators $T$ so that $[T,M_f] \in {\mathcal C}_p$ for all $f \in \Lambda(\Omega)$ with $\|Df\|_{\infty} \le 1$. The questions about ${\mathcal C}_p$-perturbations of normal operators are still covered by our results (\cite{24}, \cite{28}), however the passage of multiplication operators by Lipschitz functions would require the use of more difficult results on commutators and functional calculus, like those in \cite{2}.

\bigskip
\noindent
{\bf 6.5.} Perhaps the study of the $K$-theory of the $E\Lambda(\Omega)$ may benefit from more recent developments of bivariant $K$-theory beyond $C^*$-algebras (see \cite{14}).


\begin{thebibliography}{99}

\bibitem{1} S.~Albeverio and R.~Hoegh-Krohn, {\em Dirichlet forms and Markovian semigroups on $C^*$-algebras}, Comm. Math. Phys. {\bf 56} (1977), 173--187.

\bibitem{2} A.~B.~Aleksandrov, V.~V.~Peller, D.~S.~Potapov and F.~A.~Sukochev, {\em Functions of normal operators under perturbations}, preprint arXiv: 1008.1638 v1, 2010.

\bibitem{3} B.~Blackadar, {\em $K$-theory for operator algebras}, Math. Sci. Res. Inst. Publ. 5, Spring, 1986.

\bibitem{4} L.~G.~Brown, {\em The determinant invariant for operators with compact self-commutators}, Proc. Conf. on Operator Theory, Lecture Notes in Math., vol.~345, Springer Verlag, 1973, 210--228.

\bibitem{5} L.~G.~Brown, {\em Operator algebras and algebraic $K$-theory}, Bull. Amer. Math. Soc. (1975), 1119--1121.

\bibitem{6} L.~G.~Brown, R.~G.~Douglas and P.~A.~Fillmore, {\em Unitary equivalence modulo the compact operators and extensions of $C^*$-algebras}, Proc. Conf. on Operator Theory, Lecture Notes in Math., vol.~345, Springer Verlag, 1973, 58--127.

\bibitem{7} R.~W.~Carey and J.~D.~Pincus, {\em An exponential formula for determining functions}, Indiana Univ. Math. J. {\bf 23} (1974), 1155--1165.

\bibitem{8} R.~W.~Carey and J.~D.~Pincus, {\em Commutators, symbols and determining functions}, J. Functional Analysis {\bf 19} (1975), 50--80.

\bibitem{9} F.~Cipriani, {\em Dirichlet Forms on Non-commutative Spaces}, in U.~Franz, M.~Sch\"urmann (eds.), Quantum Potential Theory, Lecture Notes in Math. 1954, Springer Verlag 2008, 161--276.

\bibitem{10} F.~Cipriani and J.-L.~Sauvageot, {\em Derivations and square roots of Dirichlet forms}, J. Funct. Anal. {\bf 13} (2003), no.~3, 521--545.

\bibitem{11} K.~F.~Clancey, {\em Seminormal operators}, Lecture Notes in Math., vol.~742, Springer Verlag, 1979.

\bibitem{12} A.~Connes, {\em Non-commutative differential geometry}, IHES Publ. Math. {\bf 62} (1986), 257--359.

\bibitem{13} A.~Connes, {\em Non-commutative geometry}, Academic Press (1994).

\bibitem{14} J.~Cuntz and A.~Thom, {\em Algebraic $K$-theory and locally convex algebras}, Math. Ann. {\bf 334} (2006), 339--371.

\bibitem{15} J.~Duncan and S.~A.~R.~Hosseinium, {\em The second dual of a Banach algebra}, Proc. Roy. Soc. Edinburgh Sect.~A. {\bf 84} (1979), no.~3--4, 309--325.

\bibitem{16} L.~C.~Evans and R.~F.~Gariepy, {\em Measure theory and fine properties of funtions} (Studies in Advanced Mathematics) 1992, CRC Press.

\bibitem{17} J.~W.~Helton and R.~Howe, {\em Integral operators: Commutators, traces, index and homology}, Proc. Conf. on Operator Theory, Lecture Notes in Math., vol.~345, Springer Verlag, 1973, 141--209.

\bibitem{18} B.~E.~Johnson, {\em An introduction to the theory of centralizers}, Proc. London Math. Soc. {\bf 14}(3) (1964), 299--320.

\bibitem{19} O.~Lehto, {\em Univalent functions and Teichmuller spaces}, Graduate Texts in Mathematics {\bf 109}, Springer Verlag, 1986.

\bibitem{20} M.~Martin and M.~Putinar, {\em Lectures on Hyponormal Operators}, Operator Theory: Advances and Applications, vol.~39, Birkhauser, 1989.

\bibitem{21} W.~L.~Paschke, {\em $K$-theory for commutants in the Calkin algebra}, Pacific J. Math. {\bf 95}, no.~2 (1981), 427--434.

\bibitem{22} E.~M.~Stein, {\em Singular Integrals and Differentiability Properties of Functions}, Princeton University Press, 1970.

\bibitem{23} D.~Voiculescu, {\em A non-commutative Weyl--von~Neumann theorem}, Rev. Roumaine Math. Pures et Appl. {\bf 21} (1976), 97--113.

\bibitem{24} D.~V.~Voiculescu, {\em Some results on norm-ideal perturbations of Hilbert space operators I.}, J. Operator Theory {\bf 2} (1979), 3--37.

\bibitem{25} D.~V.~Voiculescu, {\em A note on quasitriangularity and trace-class self-commutators}, Acta Sci. Math. (Szeged) {\bf 42} (1980), 195--199.

\bibitem{26} D.~V.~Voiculescu, {\em Remarks on Hilbert--Schmidt perturbations of almost normal operators}, in Topics in Modern Operator Theory, Birkhauser, 1981, 311--318.

\bibitem{27} D.~V.~Voiculescu, {\em Almost normal operators modulo $\gamma_p$}, in Linear and Complex Analysis Problem Book, editors V.~P.~Havin, S.~V.~Hruscev and N.~K.~Nikolski, Lecture Notes in Math., vol.~1043, Springer Verlag, 1984, 227--230.

\bibitem{28} D.~Voiculescu, {\em Perturbations of operators, connections with singular integrals, hyperbolicity and entropy}, in Harmonic Analysis and Discrete Potential Theory (ed. M.~A.~Picardello) Plenum Press, 1992, 181--191.

\bibitem{29} N.~Weaver, {\em Lipschitz algebras}, World Scientific, 1999.

\end{thebibliography}
\end{document}